\documentclass[12pt]{article}
\usepackage{amssymb,amsmath}
\usepackage{hyperref}
\usepackage{ulem}
\usepackage{graphicx}

\begin{document}

\title{\LARGE\bf Accurate approximations for the complex error function with small imaginary argument}

%\bigskip
\author{
\normalsize\bf S. M. Abrarov\footnote{\scriptsize{Dept. Earth and Space Science and Engineering, York University, Toronto, Canada, M3J 1P3.}}\, and B. M. Quine$^{*}$\footnote{\scriptsize{Dept. Physics and Astronomy, York University, Toronto, Canada, M3J 1P3.}}}

\date{November 14, 2014}
\maketitle
%\vspace{1cm}%\bigskip

\begin{abstract}
In this paper we present two efficient approximations for the complex error function $w \left( {z} \right)$ with small imaginary argument $\operatorname{Im}{\left[ { z } \right]} <  < 1$ over the range $0 \le \operatorname{Re}{\left[ { z } \right]} \le 15$ that is commonly considered difficult for highly accurate and rapid computation. These approximations are expressed in terms of the Dawson\text{'}s integral $F\left( x \right)$ of real argument $x$ that enables their efficient implementation in a rapid algorithm. The error analysis we performed using the random input numbers $x$ and $y$ reveals that in the real and imaginary parts the average accuracy of the first approximation exceeds ${10^{ - 9}}$ and ${10^{ - 14}}$, while the average accuracy of the second approximation exceeds ${10^{ - 13}}$ and ${10^{ - 14}}$, respectively. The first approximation is slightly faster in computation. However, the second approximation provides excellent high-accuracy coverage over the required domain.
\vspace{0.25cm}
\vspace{0.25cm}
\end{abstract}

\section {Introduction}
The complex error function, also known as the Faddeeva function or the Kramp function, can be defined as \cite{Faddeeva1961, Abramowitz1972, Armstrong1972,  Poppe1990a, Poppe1990b, Schreier1992}
\begin{equation}\label{eq_1}
\begin{aligned}
w\left( z \right) &= {e^{ - {z^2}}}\left[ {1 - {\rm{erf}}\left( { - iz} \right)} \right]\\
 &= {e^{ - {z^2}}}\left( {1 + \frac{{2i}}{{\sqrt \pi  }}\int\limits_0^z {{e^{{t^2}}}dt} } \right),
\end{aligned}
\end{equation}
where $z = x + iy$ is the complex argument. Using the Fourier transforms, the real and imaginary parts of the complex error function \eqref{eq_1} can be represented as \cite{Armstrong1972, Srivastava1992, Pagnini2010}
\begin{equation}\label{eq_2}
K\left( {x,y} \right) = \frac{1}{{\sqrt \pi  }}\int\limits_0^\infty  {\exp \left( { - {t^2}/4} \right)\exp \left( { - yt} \right)\cos \left( {xt} \right)dt}, \quad \quad y > 0
\end{equation}
and
\begin{equation}\label{eq_3}
L\left( {x,y} \right) = \frac{1}{{\sqrt \pi  }}\int\limits_0^\infty  {\exp \left( { - {t^2}/4} \right)\exp \left( { - yt} \right)\sin \left( {xt} \right)dt}, \quad \quad		y > 0,
\end{equation}
respectively. The real part of the complex error function $K\left( {x,y} \right)$ is known as the Voigt function \cite{Letchworth2007}, widely used in many fields of Physics, Astronomy and Chemistry \cite{Christensen2012, Sonnenschein2012, Berk2013, Quine2013}. Combining the real \eqref{eq_2} and imaginary \eqref{eq_3} parts together as $w \left( {x,y} \right) = K \left( {x,y} \right) + iL \left( {x,y} \right)$ yields
\begin{equation}\label{eq_4}
w\left( {x,y} \right) = \frac{1}{{\sqrt \pi  }}\int\limits_0^\infty  {\exp \left( { - {t^2}/4} \right)\exp \left( { - yt} \right)\exp \left( {ixt} \right)dt}, \quad \quad y > 0.
\end{equation}
None of the integrals above can be taken analytically in closed form. Consequently, the complex error function must be computed numerically.

It should be noted that from identity \cite{McKenna1984, Zaghloul2011}
$$
w\left( { - z} \right) = 2{e^{ - {z^2}}} - w\left( z \right),
$$
it follows that only positive arguments $x$ and $y$ are sufficient in order to cover the entire complex plane. Therefore, we will imply further that both input parameters $x$ and $y$ are always positive. 

When the argument $z$ is large enough by absolute value, say $\left| {x + iy} \right| > 15$, the truncation of the Laplace continued fraction \cite{Gautschi1970, Jones1988} (see also \cite{Poppe1990a, Poppe1990b})
\begin{equation*}
w\left( z \right) = \frac{{{\mu _0}}}{{z - }}\frac{{1/2}}{{z - }}\frac{1}{{z - }}\frac{{3/2}}{{z - }}\frac{2}{{z - }}\frac{{5/2}}{{z - }} \cdots , \quad \quad {\mu _0} \equiv i/\sqrt \pi,
\end{equation*}
can be effectively used for high-accuracy and rapid computation of the complex error function $w\left( z \right)$. 

There are several approximations such as the Chiarella and Reichel approximation (equation (15) in \cite{Abrarov2014}), the Weideman's rational approximation (equation (38-I) in \cite{Weideman1994}), the exponential series approximation (equation (14) in \cite{Abrarov2011} and its modification (3) in \cite{Abrarov2012}) and the rational approximation (equation (14) in \cite{Abrarov2014}) that provide highly accurate and rapid calculation for the complex error function within domain $0 \le x \le 15$ and $10^{-6} \le y \le 15$. 

However, the highly accurate and simultaneously rapid computation of $w\left( z \right)$ at $y <  < 1$ and $\left| {x + iy} \right| \le 15$ still remains problematic \cite{Armstrong1967, Amamou2013}. Different approaches have been implemented to overcome this problem. For example, Zaghloul and Ali developed Algorithm 916 that can cover accurately this domain \cite{Zaghloul2011}. In the rapid algorithm developed by Karbach {\it{et al.}}  \cite{Karbach2014} this domain is covered by using the exponential series approximation (equation (14) in \cite{Abrarov2011}), the Taylor expansion series near singularities at $z_{n} = n \pi / \tau_{m}$ ($n$ is the index integer and $\tau_m$ is the integration cutoff) and the Laplace continued fraction. In our recent work we proposed to cover this domain by using the master-slave algorithm \cite{Abrarov2014} (see also \href{http://www.mathworks.com/matlabcentral/fileexchange/47801-the-voigt-complex-error-function}{Matlab~source~code} in \cite{MatlabCode2014}, where the master-slave approach has been implemented to cover it). In this paper we report two new approximations that also effectively resolve this problem for accurate and rapid computation.

\section{Approximations for small {\it{y}}}
A simplest way to approximate the complex error function $w \left(x,y<<1 \right)$ is to use the Maclaurin expansion series near the point $y = 0$ (see for example two approximations \eqref{eq_A.4} and \eqref{eq_A.5} in Appendix A). However, this method results in a moderate accuracy that does not exceed $10^{-6}$ in the range of a concern $\left| {x + iy} \right| \le 15$  at $y <  < 1$.

In order to obtain more efficient approximations, we may attempt to find an appropriate representation of the complex error function $w \left( {z} \right)$. Let us rewrite the equation of the complex error function \eqref{eq_4} as
\begin{equation*}
w\left( {x,y} \right) = \frac{{{e^{{y^2}}}}}{{\sqrt \pi  }}\int\limits_0^\infty  {{e^{ - {{\left( {t + 2y} \right)}^2}/4}}{e^{ixt}}dt}.
\end{equation*}
Change of the variable $u = t + 2y$ in equation above excludes the parameter $y$ from the integrand:
\begin{equation*}
\begin{aligned}
w\left( {x,y} \right) &= \frac{{{e^{{y^2}}}}}{{\sqrt \pi  }}\int\limits_{2y}^\infty  {{e^{ - {u^2}/4}}{e^{ix\left( {u - 2y} \right)}}du}
\\ &= \frac{{{e^{{y^2} - 2ixy}}}}{{\sqrt \pi  }}\int\limits_{2y}^\infty  {{e^{ - {u^2}/4}}{e^{ixu}}du.} 
\end{aligned}
\end{equation*}
Consequently, this equation can be rearranged as
\begin{equation}\label{eq_5}
w\left( {x,y} \right) = \frac{{{e^{{y^2} - 2ixy}}}}{{\sqrt \pi  }}\left( {\int\limits_0^\infty  {{e^{ - {u^2}/4}}{e^{ixu}}du} - \int\limits_0^{2y} {{e^{ - {u^2}/4}}{e^{ixu}}du}} \right).
\end{equation}

The first integral in the equation above is
\begin{equation*}
\int\limits_0^\infty  {{e^{ - {u^2}/4}}{e^{ixu}}du}  = {e^{ - {x^2}}}\sqrt \pi + 2i\,F\left( x \right),
\end{equation*}
where
$$
F\left( x \right) = \frac{1}{2}\int\limits_0^\infty  {{e^{ - {u^2}/4}}\sin \left( {xu} \right)du}
$$
is the Dawson\text{'}s integral. As the argument $x$ in $F \left( x \right)$ is real, its computation is not difficult and several efficient approximations that can provide rapid and highly accurate computation are reported in literature \cite{Rybicki1989, McCabe1974, Cody1970}. Furthermore, the latest versions of Matlab support the built-in function {\sffamily{dawson(x)}} that computes very rapidly the Dawson\text{'}s integral of real argument.

The advantage of the equation \eqref{eq_5} is that the integration range in the second integral at $y << 1$ is very narrow. Consequently, the second integral in equation \eqref{eq_5} makes only a minor contribution to the complex error function $w\left({z}\right)$. As a result, this minimizes the error in computation.

The representation of the complex error function in form of \eqref{eq_5} provides enormous flexibility to approximate it. Here we show two efficient derivations that directly follow from equation \eqref{eq_5}.

As the upper limit of the second integral of equation \eqref{eq_5} is finite,  its evaluation is not straightforward (see Appendix B). However, this problem can be overcome by Maclaurin expansion of the exponential function:
$$
{e^{ - {u^2}/4}} = 1 - \frac{{{u^2}}}{4} + \frac{{{u^4}}}{{32}} - \frac{{{u^6}}}{{384}} + O\left( {{u^8}} \right).
$$

The first approximation is obtained by substituting only the first term of this expansion into the second integral of equation \eqref{eq_5} that yields
\begin{equation}\label{eq_6}
w\left( {x,y <  < 1} \right) \approx {e^{{{\left( {ix - y} \right)}^2}}}\left[ {1 + \frac{{i{e^{{x^2}}}}}{{\sqrt \pi  }}\left( {2F\left( x \right) - \frac{{1 - {e^{2ixy}}}}{x}} \right)} \right].
\end{equation}
We will refer to this equation as the basic approximation. This equation is highly accurate while $x \lesssim 3$.

The second approximation is obtained by substituting the first two terms of the expansion above into second integral of the equation \eqref{eq_5}:
\small
\begin{equation}\label{eq_7}
\begin{aligned}
w\left( {x,y <  < 1} \right) &\approx {e^{{{\left( {ix - y} \right)}^2}}} \left[ 1 + \frac{{i{e^{{x^2}}}}}{{\sqrt \pi  }} \right.
\\ &\times \left. \left( 2F\left( x \right) - \frac{{{{x^2} + 0.5 + e^{2ixy}}\left( {{x^2}\left( {{y^2} - 1} \right) + ixy - 0.5} \right)}}{{{x^3}}}  \right) \right].
\end{aligned}
\end{equation}
\normalsize
Further, we will refer to this equation as the main approximation.

Although the basic approximation \eqref{eq_6} is slightly faster in computation, the main approximation \eqref{eq_7} provides essentially higher accuracy as it will be shown in the section $4$.

\section{Algorithmic implementation}
Approximations \eqref{eq_6} and \eqref{eq_7} cannot be employed directly in the computation flow since at $x \to 0$ the denominators $x$ in approximation \eqref{eq_6} and $x^{3}$ in approximation \eqref{eq_7} blow up the ratios resulting to computational overflow. However, this problem is very easy to resolve by using, for example, the following supplement:
\begin{equation}\label{eq_8}
w\left( {x <  < 1,y <  < 1} \right) \approx \left( {1 - {\left( {x + iy} \right)^2}} \right)\left( {1 + \frac{{2i \left( {x + iy} \right)}}{{\sqrt \pi  }}} \right).
\end{equation}
In particular, the computation flow for the basic and the main approximations can be maintained as
\begin{equation*}
w\left( {x,y <  < 1} \right) \approx \left\{ 
\begin{aligned}
{\rm{Eq.(6)}}\, {\rm{or}}  \, {\rm{Eq.(7),}} &\quad {10^{-4}} < x \le 15 \cap y \le {10^{ - 6}}\\
{\rm{Eq.(8)}}, &\quad 0 \le x \le{10^{-4}} \cap y \le {10^{ - 6}}.
\end{aligned} \right.
\end{equation*}
Thus, according to this scheme if  $x > {10^{-4}}$ the complex error function is computed either by equation \eqref{eq_6} or by equation \eqref{eq_7}. Otherwise, if $x \le {10^{-4}}$ it is computed by equation \eqref{eq_8}.

The approximation \eqref{eq_8} had been used already in our recently published \href{http://www.mathworks.com/matlabcentral/fileexchange/47801-the-voigt-complex-error-function}{Matlab~source~code} \cite{MatlabCode2014} to resolve a similar problem that occurs at $\left| {x + iy} \right| \to 0$.  The derivation of the approximation \eqref{eq_8} is shown in Appendix A.

\section{Error analysis}
Define the relative errors for the real and imaginary parts of the complex error function in forms
$$
{\Delta _{{\mathop{\rm Re}\nolimits} }} = \left| {\frac{{{\mathop{\rm Re}\nolimits} \left[ {w\left( z \right)} \right] - {\mathop{\rm Re}\nolimits} \left[ {{w_{ref.}}\left( z \right)} \right]}}{{{\mathop{\rm Re}\nolimits} \left[ {{w_{ref.}}\left( z \right)} \right]}}} \right|
$$
and
$$
{\Delta _{{\mathop{\rm Im}\nolimits} }} = \left| {\frac{{{\mathop{\rm Im}\nolimits} \left[ {w\left( z \right)} \right] - {\mathop{\rm Im}\nolimits} \left[ {{w_{ref.}}\left( z \right)} \right]}}{{{\mathop{\rm Im}\nolimits} \left[ {{w_{ref.}}\left( z \right)} \right]}}} \right|,
$$
respectively, where ${w_{ref.}}\left( z \right)$ is the reference. The highly accurate reference values can be obtained according to equation \eqref{eq_1} by using the latest versions of Wolfram Mathematica that supports error function of complex argument.

Figure 1a shows the logarithm of relative error for the real part of the basic approximation \eqref{eq_6} in the domain $10^{-4} \le x \le 15$ and $y \le 10^{-6}$. As we can see, at $x \lesssim 3$ the approximation \eqref{eq_6} is highly accurate and provides accuracy better than $10^{-12}$. Although the accuracy of the approximation \eqref{eq_6} deteriorates as $x$ increases, it, nevertheless, still remains high and better than $10^{-9}$.

\newpage
\begin{figure}[ht]
\begin{center}
\includegraphics[width=17pc]{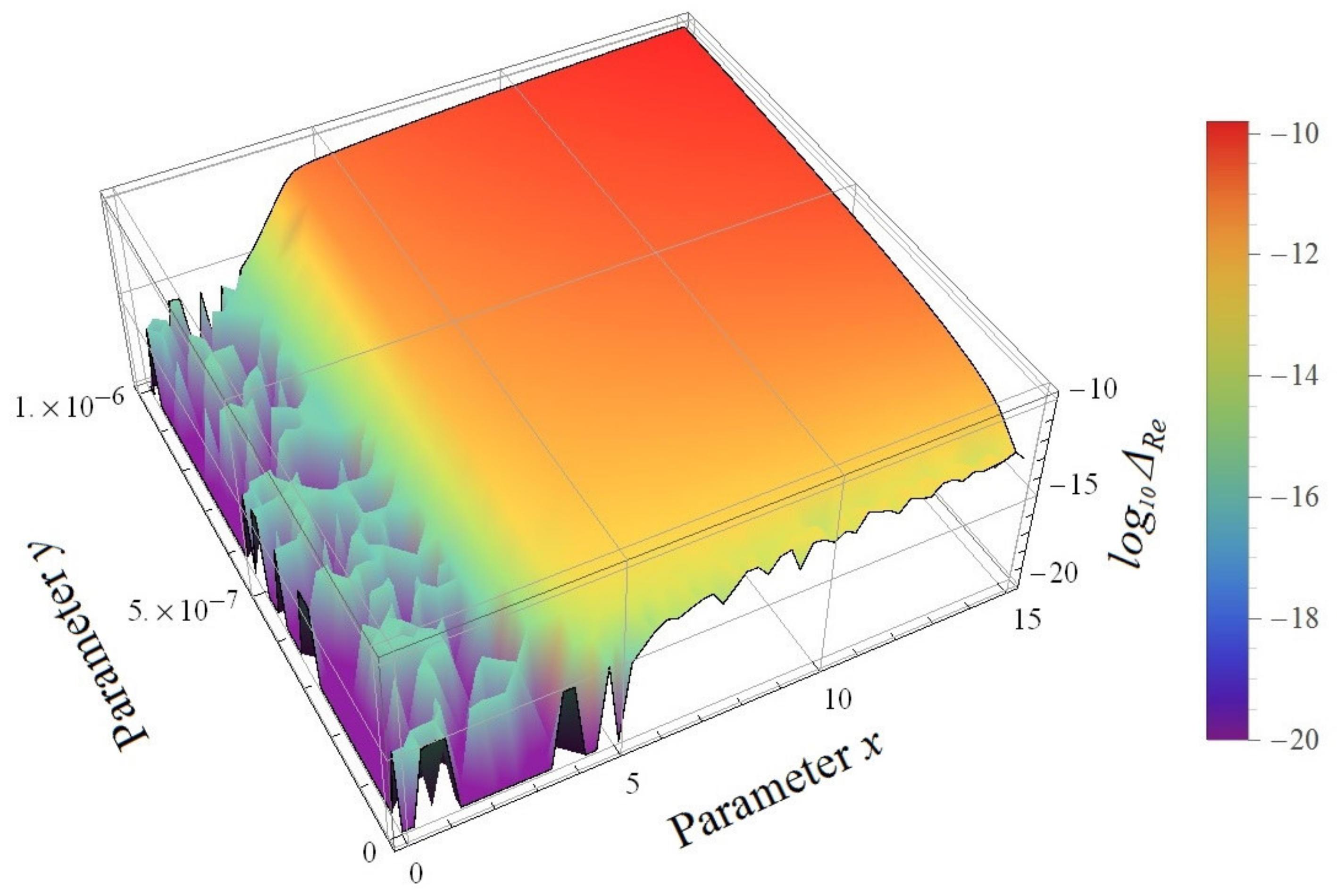}\hspace{2pc}%
\begin{minipage}[b]{24pc}
\vspace{0.3cm}
{\sffamily {\bf{Fig. 1a.}} Logarithm of relative error $log_{10}\Delta_{\operatorname{Re}}$ for the real part of the basic approximation \eqref{eq_6} over the domain $10^{-4} \le x \le 15$ and $y \le 10^{-6}$.}
\end{minipage}
\end{center}
\end{figure}

\begin{figure}[ht]
\begin{center}
\includegraphics[width=17pc]{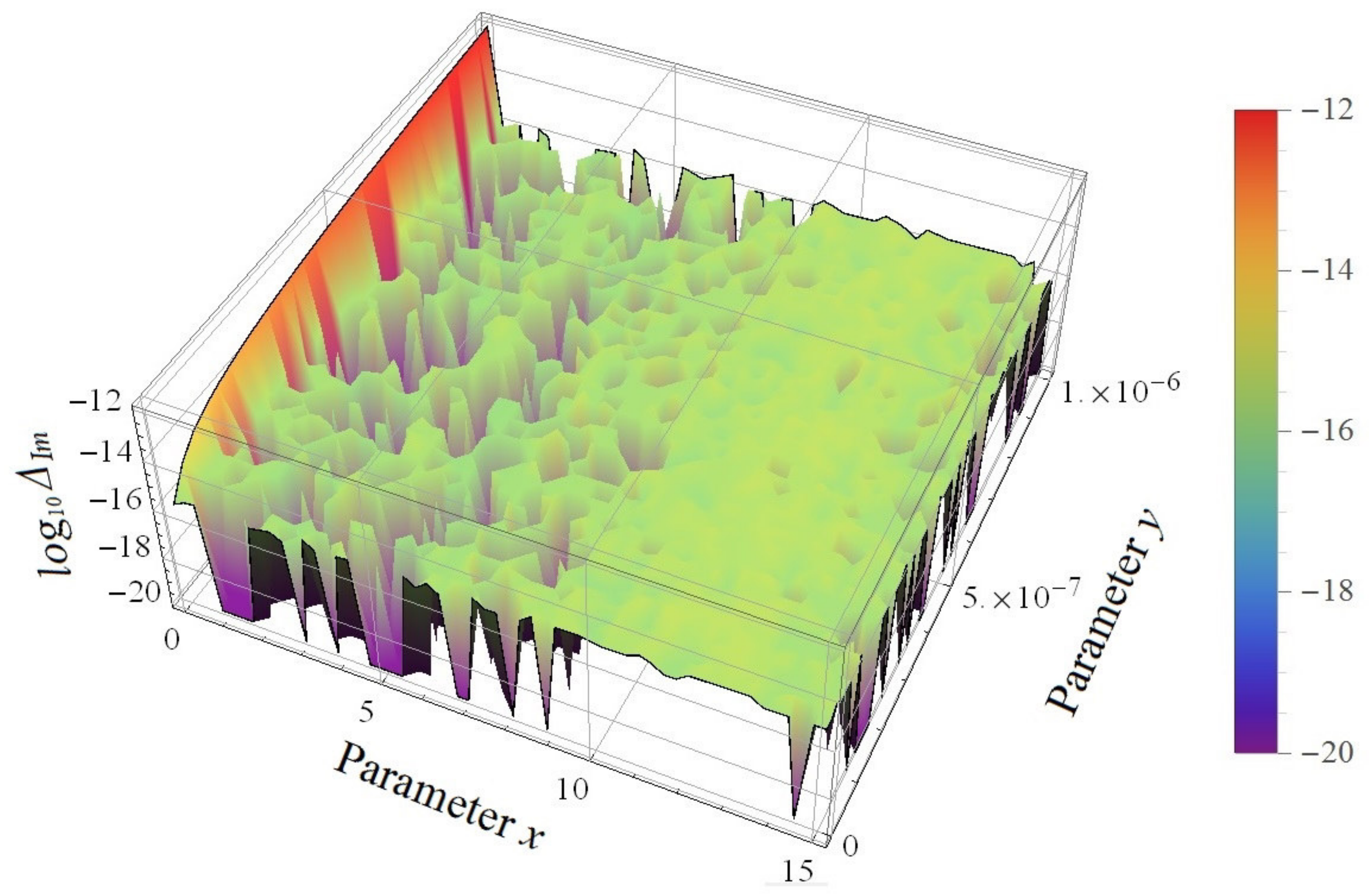}\hspace{2pc}%
\begin{minipage}[b]{24pc}
\vspace{0.3cm}
{\sffamily {\bf{Fig. 1b.}} Logarithm of relative error $log_{10}\Delta_{\operatorname{Im}}$ for the imaginary part of the basic approximation \eqref{eq_6} over the domain $10^{-4} \le x \le 15$ and $y \le 10^{-6}$.}
\end{minipage}
\end{center}
\end{figure}

Figure 1b illustrates the logarithm of relative error for the imaginary part of the basic approximation \eqref{eq_6} in the domain $10^{-4} \le x \le 15$ and $y \le 10^{-6}$. As one can see, the imaginary part of the approximation \eqref{eq_6} is highly accurate. Specifically, the accuracy is better than $10^{-14}$ almost over all of the domain. However, at small $x \lesssim 10^{-3}$ the accuracy deteriorates as it can be seen by red color along the left edge of $y$-axis.

\newpage
\begin{figure}[ht]
\begin{center}
\includegraphics[width=18pc]{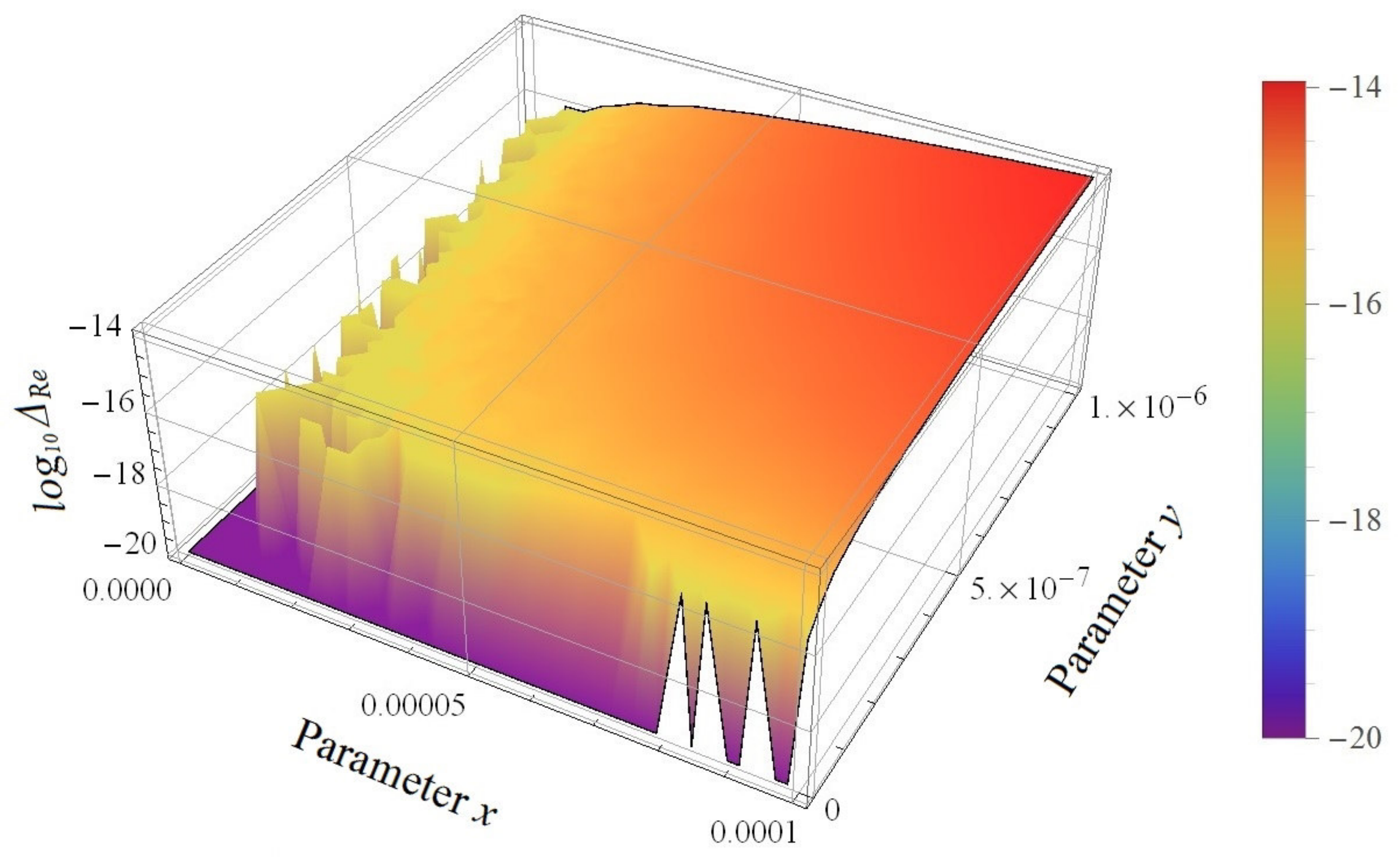}\hspace{2pc}%
\begin{minipage}[b]{24pc}
\vspace{0.3cm}
{\sffamily {\bf{Fig. 2a.}} Logarithm of relative error $log_{10}\Delta_{\operatorname{Re}}$ for the real part of the supplementary approximation \eqref{eq_8} over the small domain $0 \le x \le 10^{-4}$ and $y \le 10^{-6}$.}
\end{minipage}
\end{center}
\end{figure}

\begin{figure}[ht]
\begin{center}
\includegraphics[width=18pc]{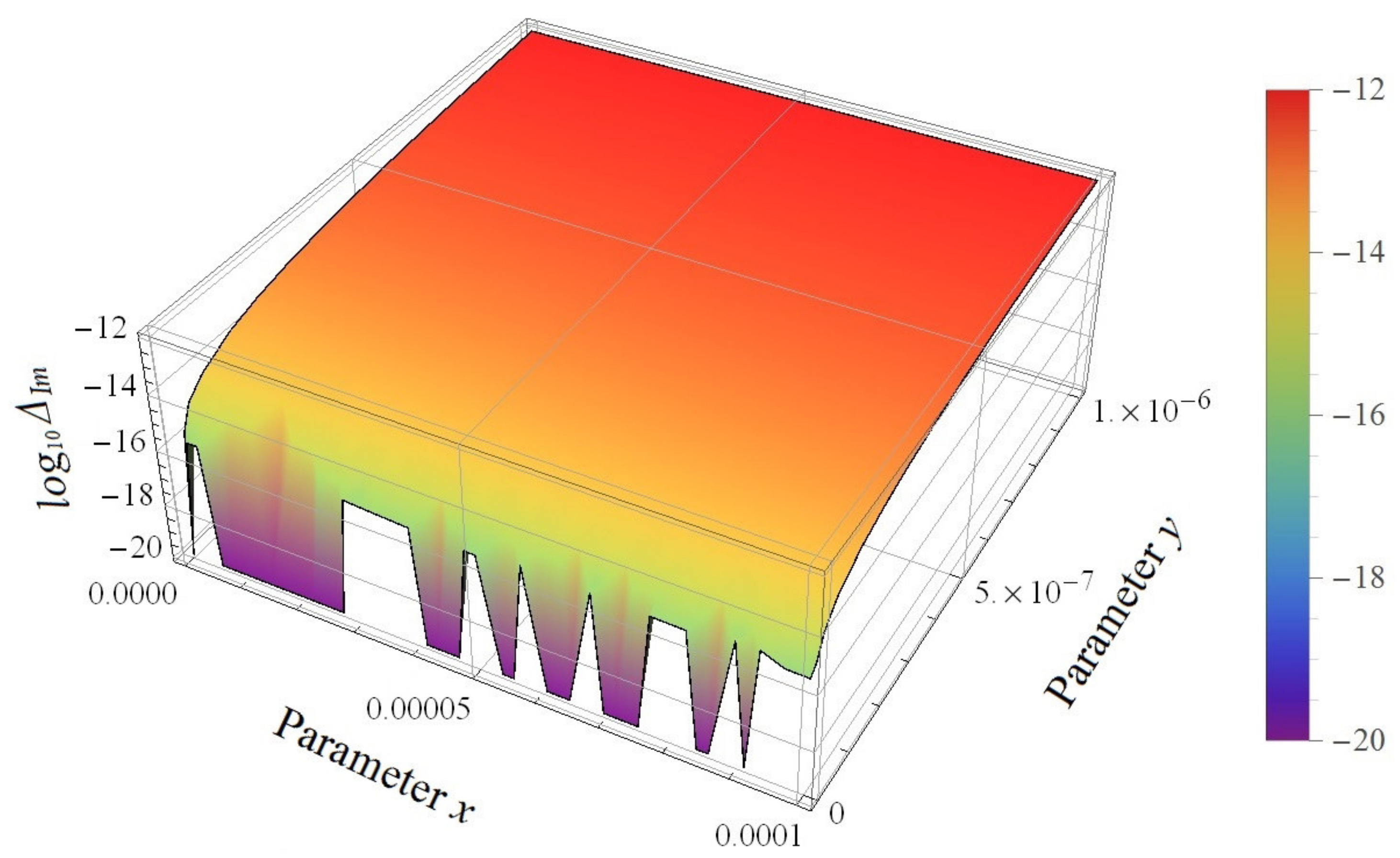}\hspace{2pc}%
\begin{minipage}[b]{24pc}
\vspace{0.3cm}
{\sffamily {\bf{Fig. 2b.}} Logarithm of relative error $log_{10}\Delta_{\operatorname{Im}}$ for the imaginary part of the supplementary approximation \eqref{eq_8} over the small domain $0 \le x \le 10^{-4}$ and $y \le 10^{-6}$.}
\end{minipage}
\end{center}
\end{figure}

Figure 2a illustrates the logarithm of relative error $log_{10}\Delta_{\operatorname{Re}}$ for the real part of the supplementary approximation \eqref{eq_8} over the small domain $0 \le x \le 10^{-4}$ and $y \le 10^{-6}$. One can see that the real part over this domain is highly accurate and better than $10^{-13}$.

Figure 2b depicts the logarithm of relative error $log_{10}\Delta_{\operatorname{Im}}$ for the imaginary part of the supplementary approximation \eqref{eq_8} over the small domain $0 \le x \le 10^{-4}$ and $y \le 10^{-6}$. The accuracy of the imaginary part is also high and better than $10^{-12}$.

From Figs. 2a and 2b we can conclude that the basic approximation \eqref{eq_6} alone may be used practically as the accuracy better than $10^{-9}$ is more than enough for the most applications.

%\newpage
\begin{figure}[ht]
\begin{center}
\includegraphics[width=15pc]{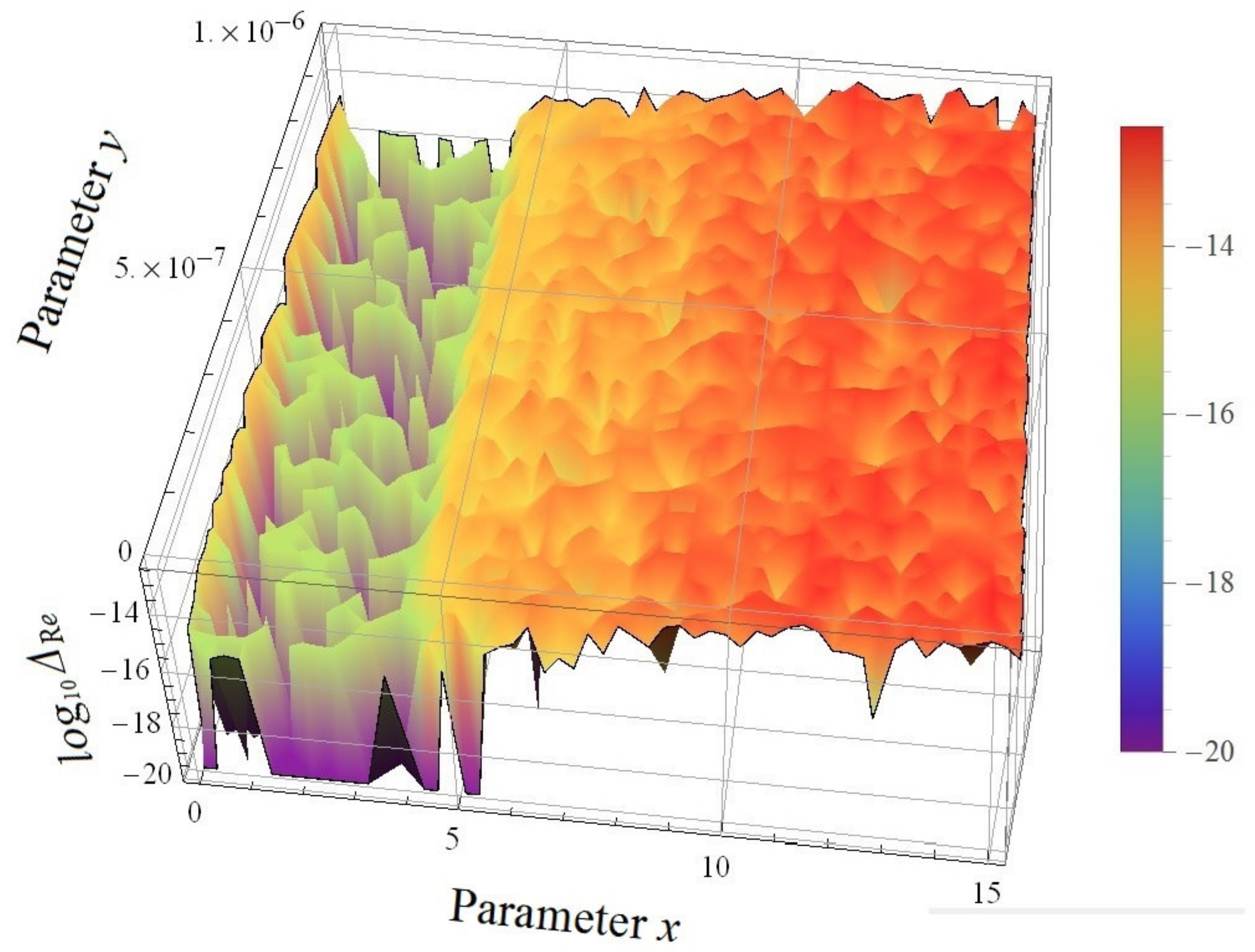}\hspace{2pc}%
\begin{minipage}[b]{24pc}
\vspace{0.3cm}
{\sffamily {\bf{Fig. 3a.}} Logarithm of relative error $log_{10}\Delta_{\operatorname{Re}}$ for the real part of the main approximation \eqref{eq_7} over the domain $10^{-4} \le x \le 15$ and $y \le 10^{-6}$.}
\end{minipage}
\end{center}
\end{figure}

\begin{figure}[ht]
\begin{center}
\includegraphics[width=18pc]{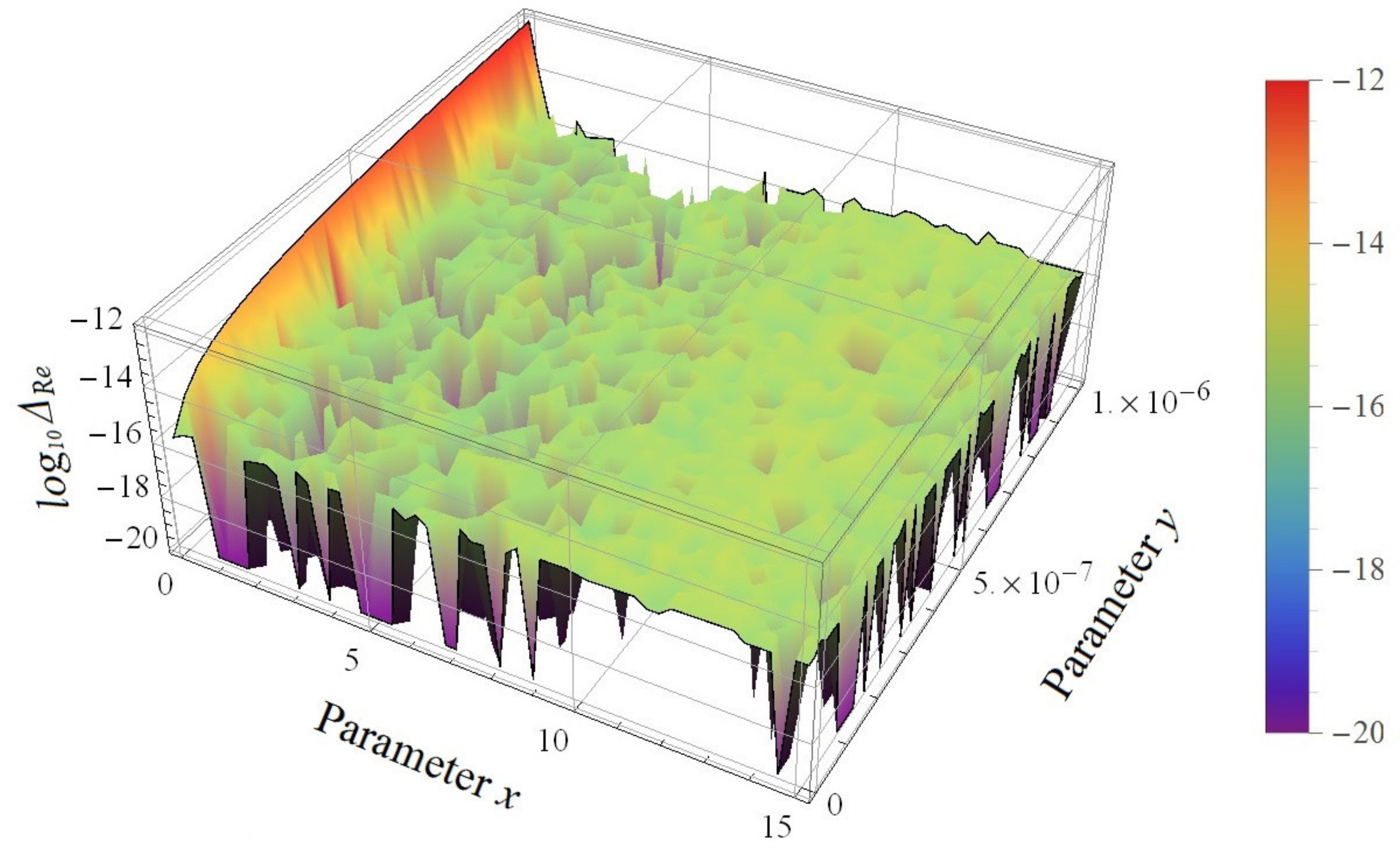}\hspace{2pc}%
\begin{minipage}[b]{24pc}
\vspace{0.3cm}
{\sffamily {\bf{Fig. 3b.}} Logarithm of relative error $log_{10}\Delta_{\operatorname{Im}}$ for the imaginary part of the main approximation \eqref{eq_7} over the domain $10^{-4} \le x \le 15$ and $y \le 10^{-6}$.}
\end{minipage}
\end{center}
\end{figure}

The accuracy of the complex error function $w \left( {x,y<<1} \right)$ can be further improved by using the main approximation \eqref{eq_7}. Figure 3a depicts the logarithm of relative error for the real part of the approximation \eqref{eq_7} in the domain $10^{-4} \le x \le 15$ and $y \le 10^{-6}$. As it can be seen from this figure, in the real part the accuracy of the approximation \eqref{eq_7} is better than $10^{-13}$.

Figure 3b shows the logarithm of relative error for the imaginary part of the main approximation \eqref{eq_7} in the domain $10^{-4} \le x \le 15$ and $y \le 10^{-6}$. Comparing Figs. 3b with 1b one can observe that the graphs are nearly same. This signifies that the behavior of the basic and the main approximations resemble in their imaginary parts. Similar to Fig. 1b, the accuracy in the imaginary part is better than $10^{-14}$ almost over all domain.

To determine the average accuracy for both approximations, we performed calculations with randomly chosen input variables $x$ and $y$. It has been found that in the real and imaginary parts the basic approximation \eqref{eq_6} provides average accuracy exceeding $10^{-9}$ and $10^{-14}$, while the main approximation \eqref{eq_7} provides average accuracy exceeding $10^{-13}$ and $10^{-14}$, respectively.

\section{Conclusion}
We present two efficient approximations for the complex error function that can be used for computation at small $y <  < 1$. As these approximation are expressed in terms of the Dawson\text{'}s integral $F\left( x \right)$ of real argument $x$, they are rapid in computation. Error analysis performed with randomly taken input variables $x$ and $y$ reveals that in the real and imaginary parts the basic approximation \eqref{eq_6} provides average accuracy exceeding ${10^{ - 9}}$ and ${10^{ - 14}}$, while the main approximation \eqref{eq_7} provides average accuracy exceeding ${10^{ - 13}}$ and ${10^{ - 14}}$, respectively. Although the basic equation \eqref{eq_6} is slightly faster in computation, the main approximation \eqref{eq_7} provides excellent high-accuracy coverage.

\section*{Appendix A}
Consider the following Maclaurin expansions at $y = 0$:
\[\label{eq_A.1}
\tag{A.1}
{e^{ - {{\left( {x + iy} \right)}^2}}} = {e^{ - {x^2}}} - 2ixy{e^{ - {x^2}}} + O\left( {{y^2}} \right),
\]
\[\label{eq_A.2}
\tag{A.2}
{\rm{erf}}\left( { - i\left( {x + iy} \right)} \right) =  - i{\rm{erfi}}\left( x \right) + \frac{{2y}}{{\sqrt \pi  }}{e^{x^2}} + O\left( {{y^2}} \right)
\]
and
\footnotesize
\[\label{eq_A.3}
\tag{A.3}
\begin{aligned}
w\left( {x,y} \right) &\equiv {e^{ - {{\left( {x + iy} \right)}^2}}}\left[ {1 - {\rm{erf}}\left( { - i\left( {x + iy} \right)} \right)} \right] 
\\&= {e^{ - {x^2}}}\left( {1 + {i\rm{erfi}}\left( {x} \right)} \right) + 2y\left( {x{e^{ - {x^2}}}{\rm{erfi}}\left( x \right) - ix{e^{ - {x^2}}} - \frac{1}{{\sqrt \pi  }}} \right) + O\left( {{y^2}} \right),
\end{aligned}
\]
\normalsize
Substituting \eqref{eq_A.1} and \eqref{eq_A.2} into identity \eqref{eq_1} we obtain
\[\label{eq_A.4}
\tag{A.4}
w\left( {x,y <  < 1} \right) \approx \left( {1 - 2ixy} \right)\left( {{e^{ - {x^2}}}\left( {1 + {\rm{erf}}\left( {ix} \right)} \right) - \frac{{2y}}{{\sqrt \pi  }}} \right),
\]
while equation \eqref{eq_A.3} immediately results in
\[\label{eq_A.5}
\tag{A.5}
w\left( {x,y <  < 1} \right) \approx {e^{ - {x^2}}}\left( {1 + {\rm{erf}}\left( {ix} \right)} \right)\left( {1 - 2ixy} \right) - \frac{{2y}}{{\sqrt \pi  }}.
\]

The approximations \eqref{eq_A.4} and \eqref{eq_A.5} are rapid in computation because the error function of imaginary argument ${\rm{erf}}\left( {ix} \right)$ can be expressed in terms of the Dawson's integral of real argument $F \left( x \right)$ as follows
\[
{\rm{erf}}\left( {ix} \right) \equiv {\frac{2i}{\sqrt{\pi}}} e^{x^{2}} F \left( x \right).
\]
However, the approximations \eqref{eq_A.4} and \eqref{eq_A.5} are moderately accurate and should be used when the requirement for accuracy does not exceed $10^{-6}$.

Same technique can be used to obtain approximation for $w\left( {x < < 1, y < < 1} \right)$ from the expansions near the point $z=0$:
\[
\exp \left( { - {z^2}} \right) = 1 - {z^2} + O\left( {{z^4}} \right),
\]
and
\[
{\rm{erf}}\left( { - iz} \right) =  - \frac{{2iz}}{{\sqrt \pi  }} + O\left( {{z^3}} \right),
\]
In particular, substituting these expansions into identity \eqref{eq_1} leads to the approximation \eqref{eq_8}.

\section*{Appendix B}
The second integral of equation \eqref{eq_5} can be expressed as
\[
\int\limits_0^{2y} {{e^{ - {u^2}/4}}{e^{ixu}}du}  = \sqrt \pi  {e^{ - {x^2}}}\left[ {{\rm{erf}}\left( {ix} \right) - {\rm{erf}}\left( {ix - y} \right)} \right].
\]
The argument in the first error function ${\rm{erf}}\left({ix}\right)$ is purely imaginary and, therefore, it can be expressed in term of the Dawson\text{'}s integral of real argument. However, since the argument in the second error function ${\rm{erf}}\left({ix - y}\right)$ is complex, it needs a numerical solution (see for example equations (10) and (11) in \cite{Zaghloul2011}). When the high accuracy is not required, the rapid approximation for the error function can obtained, for example, by rearranging the expansion series \eqref{eq_A.2} as
\[
\begin{aligned}
{\rm{erf}}\left( {ix - y} \right) &\approx i{\rm{erfi}}\left( x \right) - \frac{{2y}}{{\sqrt \pi  }}{e^{{x^2}}}
\\ &\equiv {\rm{erf}}\left( {ix} \right) - \frac{{2y}}{{\sqrt \pi  }}{e^{{x^2}}}, \quad\quad y << 1.
\end{aligned}
\]

We can also approximate this integral by taking into account that $y << 1$. In particular, for small interval we can write $du \approx 2y$ and take the mid-point in that interval equal to $y$:
$$
\int\limits_0^{2y} {{e^{ - {u^2}/4}}{e^{ixu}}du}  \approx 2y{e^{ - {y^2}/4}}{e^{ixy}}.
$$
Alternatively, according to the trapezoidal rule we get
$$
\int\limits_0^{2y} {{e^{ - {u^2}/4}}{e^{ixu}}du}  \approx y\left( {1 + {e^{ - {y^2}}}{e^{2ixy}}} \right).
$$
Once again, these approximations should be applied only when the high-accuracy is not a concern in computation.

The second integral from equation \eqref{eq_5} can be approximated more precisely by taking integration by parts
\[
\label{eq_B.1}
\tag{B.1}
\int\limits_0^{2y} {{e^{ - {u^2 / 4}}}{e^{ixu}}} du = \left. {\sqrt \pi \, {\rm{erf}}\left( {\frac{u}{2}} \right){e^{ixu}}} \right|_0^{2y} - \sqrt \pi  ix\int\limits_0^{2y} {{\rm{erf}}\left( {\frac{u}{2}} \right){e^{ixt}}} du
\]
and then substituting the first few terms from the following expansion
\[
{\rm{erf}}\left( {\frac{u}{2}} \right) = \frac{u}{{\sqrt \pi  }} - \frac{{{u^3}}}{{12\sqrt \pi  }} + \frac{{{u^5}}}{{160\sqrt \pi  }} - \frac{{{u^7}}}{{2688\sqrt \pi  }} + O\left( {{u^9}} \right)
\]
into the right integral of equation \eqref{eq_B.1}. We will consider the simplest case when only the first term from this expansion is substituted. Following this procedure and using equation \eqref{eq_5} we have
\small
\[
\label{eq_B.2}
\tag{B.2}
\begin{aligned}
w\left( {x,y <  < 1} \right) &\approx {e^{{{\left( {ix - y} \right)}^2}}} \left[ 1 + \frac{{i{e^{{x^2}}}}}{{\sqrt \pi  }} \right.
\\ &\times \left. \left( {2F\left( x \right) - \frac{{1 - {e^{2ixy}}\left( {1 + ix \left( \sqrt \pi \, {\rm{erf}}\left( y \right) - 2y \right)} \right)}}{x}} \right) \right].
\end{aligned}
\]
\normalsize
The consistency between approximations \eqref{eq_6} and \eqref{eq_B.2} becomes evident from the fact that ${\rm{erf}} \left( {y << 1} \right) \approx 2y / \sqrt \pi$. Although approximation \eqref{eq_B.2} is more accurate than the basic approximation \eqref{eq_6}, it is slightly slower in computation at large size of the input array due to presence of the error function ${\rm{erf}} \left( {y} \right)$. Since the approximation \eqref{eq_B.2} provides nearly same accuracy as the main approximation \eqref{eq_7}, both of them can be used when the high-accuracy computation is required.

\small
\section*{Acknowledgments}
This work is supported by National Research Council of Canada, Thoth Technology Inc. and York University. The authors wish to thank Prof. Ian McDade and Dr. Brian Solheim for discussions and constructive suggestions.

%\bigskip
%\newpage


\small
\begin{thebibliography}{9}

\bibitem{Faddeeva1961}
V. N. Faddeyeva and N. M. Terent\text{’}ev, Tables of the probability integral for complex argument, Pergamon Press, Oxford, 1961.

\bibitem{Abramowitz1972}
M. Abramowitz, and I. A. Stegun, Error function and Fresnel integrals. Handbook of mathematical functions with formulas, graphs, and mathematical tables. $9^{th}$ Ed. Dover, New York, 1972, 297-309.

\bibitem{Armstrong1972}
B. H. Armstrong and B. W. Nicholls, Emission, absorption and transfer of radiation in heated atmospheres, Pergamon Press, New York, 1972.

\bibitem{Poppe1990a}
G. P. M. Poppe and C. M. J. Wijers, More efficient computation of the complex error function, \href{http://doi.org/10.1145/77626.77629}{ACM Transact. Math. Software, 16 (1990) 38-46}.

\bibitem{Poppe1990b}
G. P. M. Poppe and C. M. J. Wijers, Algorithm 680: evaluation of the complex error function, \href{http://doi.org/10.1145/77626.77630}{ACM Transact. Math. Software, 16 (1990) 47}.

\bibitem{Schreier1992}
F. Schreier, The Voigt and complex error function: A comparison of computational methods, \href{http://doi.org/10.1016/0022-4073(92)90139-U}{J. Quant. Spectrosc. Radiat. Transfer, 48 (1992) 743-762}.

\bibitem{Srivastava1992}
H. M. Srivastava and M. P. Chen, Some unified presentations of the Voigt functions, \href{http://doi.org/10.1007/BF00653260}{Astrophys. Space Sci., 192 (1992) 63-74}.

\bibitem{Pagnini2010}
G. Pagnini and F. Mainardi, Evolution equations for the probabilistic generalization of the Voigt profile function, \href{http://dx.doi.org/10.1016/j.cam.2008.04.040}{J. Comput. Appl. Math., 233 (2010) 1590-1595}.

\bibitem{Letchworth2007}
K. L. Letchworth and D. C. Benner, Rapid and accurate calculation of the Voigt function, \href{http://dx.doi.org/10.1016/j.jqsrt.2007.01.052}{J. Quantit. Spectrosc. Radiat. Transfer, 107 (2007) 173-192}.

\bibitem{Christensen2012}
L. E. Christensen, G. D. Spiers, R. T. Menzies and J. C. Jacob, Tunable laser spectroscopy of $\rm{CO}_{2}$ near $2.05 \, {\mu}m$: Atmospheric retrieval biases due to neglecting line-mixing, \href{http://doi.org/10.1016/j.jqsrt.2012.02.031}{J. Quantit. Spectrosc. Radiat. Transfer, 113 (2012) 739-748}.

\bibitem{Sonnenschein2012}
V. Sonnenschein, S. Raeder, A. Hakimi, I. D. Moore and K. Wendt, Determination of the ground-state hyperfine structure in neutral $^{229}{\rm{Th}}$, \href{http://dx.doi.org/10.1088/0953-4075/45/16/165005}{J. Phys. B: At. Mol. Opt. Phys. 45 (2012) 165005}.

\bibitem{Berk2013}
A. Berk, Voigt equivalent widths and spectral-bin single-line transmittances: Exact expansions and the MODTRAN{\circledR}5 implementation, \href{http://dx.doi.org/10.1016/j.jqsrt.2012.11.026}{J. Quantit. Spectrosc. Radiat. Transfer, 118 (2013) 102-120}.

\bibitem{Quine2013}
B. M. Quine and S. M. Abrarov, Application of the spectrally integrated Voigt function to line-by-line radiative transfer modelling,
\href{http://doi.org/10.1016/j.jqsrt.2013.04.020}{J. Quantit. Spectrosc. Radiat. Transfer, 127 (2013) 37-48}.

\bibitem{McKenna1984}
S. J. McKenna, A method of computing the complex probability function and other related functions over the whole complex plane, \href{http://doi.org/10.1007/BF00649615}{Astrophys. Space Sci., 107 (1984) 71-83}.

\bibitem{Zaghloul2011}
M. R. Zaghloul and A. N. Ali, Algorithm 916: computing the Faddeyeva and Voigt functions, \href{http://doi.org/10.1145/2049673.2049679}{ACM Trans. Math. Software 38 (2011) 15:1-15:22}.

\bibitem{Gautschi1970}
W. Gautschi, Efficient computation of the complex error function, \href{http://dx.doi.org/10.1137/0707012}{SIAM J. Numer. Anal., 7 (1970) 187-198}.

\bibitem{Jones1988}
W. B. Jones and W. J. Thron, Continued fractions in numerical analysis, \href{http://doi.org/10.1016/0168-9274(83)90002-8}{Appl. Num. Math., 4 (1988) 143-230}.

\bibitem{Abrarov2014}
S. M. Abrarov and B. M. Quine, Master-slave algorithm for highly accurate and rapid computation of the Voigt/complex error function, \href{http://doi.org/10.5539/jmr.v6n2p104}{J.~Math. Research, 6 (2014) 104-119}.

\bibitem{Weideman1994}
J. A. C. Weideman, Computation of the complex error function, \href{http://dx.doi.org/10.1137/0731077}{SIAM J. Numer. Anal., 31 (1994) 1497-1518}.

\bibitem{Abrarov2011}
S. M. Abrarov and B. M. Quine, Efficient algorithmic implementation of the Voigt/complex error function based on exponential series approximation, \href{http://doi.org/10.1016/j.amc.2011.06.072}{Appl. Math. Comput. 218 (2011) 1894-1902}.

\bibitem{Abrarov2012}
S. M. Abrarov and B. M. Quine,  On the Fourier expansion method for highly accurate computation of the Voigt/complex error function in a rapid algorithm, \href{http://arxiv.org/pdf/1205.1768v1.pdf}{arXiv:1205.1768}.

\bibitem{Karbach2014}
T. M. Karbach, G. Raven and M. Schiller, Decay time integrals in neutral meson mixing and their efficient evaluation, \href{http://arxiv.org/pdf/1407.0748v1.pdf}{arXiv:1407.0748}.

\bibitem{Armstrong1967}
B. H. Armstrong, Spectrum line profiles: the Voigt function, \href{http://doi.org/10.1016/0022-4073(67)90057-X}{J. Quant. Spectrosc. Radiat. Transfer. 7 (1967) 61-88}.

\bibitem{Amamou2013}
H. Amamou, B. Ferhat and A. Bois, Calculation of the Voigt function in the region of very small values of the parameter a where the calculation is notoriously difficult, \href{10.4236/ajac.2013.412087}{Amer. J.  Anal. Chem., 4 (2013) 725-731}.

\bibitem{MatlabCode2014}
The Matlab source code for computation of the Voigt/complex error function can be downloaded here:
\href{http://www.mathworks.com/matlabcentral/fileexchange/47801-the-voigt-complex-error-function}{Matlab Central, file ID: \#47801, submitted~on Sept. 10, 2014}.

\bibitem{Rybicki1989}
G. B. Rybicki, Dawson\text{'}s integral and the sampling theorem, \href{http://dx.doi.org/10.1063/1.4822832}{Comp. Phys., 3 (1989) 85-87}.

\bibitem{McCabe1974}
J. H. McCabe, A continued fraction expansion with a truncation error estimate for Dawson\text{'}s integral, \href{http://dx.doi.org/10.1090/S0025-5718-1974-0371020-3#sthash.Ou4P83VO.dpuf}{Math. Comp. 127 (1974) 811-816}.

\bibitem{Cody1970}
W. J. Cody, K. A. Paciorek and H. C. Thacher, Chebyshev approximations for Dawson\text{'}s integral, \href{http://dx.doi.org/10.1090/S0025-5718-1970-0258236-8}{Math. Comp. 24 (1970) 171-178}.

\end{thebibliography}
\end{document}